\documentclass[12pt,reqno]{amsart}
\usepackage{enumerate, latexsym, amsmath, amsfonts, amssymb, amsthm, graphicx, color}
\def\pmod #1{\ ({\rm{mod}}\ #1)}
\def\Z{\Bbb Z}
\def\N{\Bbb N}

\def\l{\left}
\def\r{\right}
\def\bg{\bigg}
\def\({\bg(}
\def\){\bg)}
\def\t{\text}
\def\f{\frac}

\def\ls{\leqslant}
\def\gs{\geqslant}
\def\se {\subseteq}
\def\sm{\setminus}

\def\eq{\equiv}

\def\da{\delta}

\def\Proof{\noindent{\it Proof}}

\theoremstyle{plain}
\newtheorem{theorem}{Theorem}

\newtheorem{lemma}{Lemma}
\newtheorem{corollary}{Corollary}

\theoremstyle{definition}

\theoremstyle{remark}
\newtheorem{remark}{Remark}

\makeatletter
\@namedef{subjclassname@2010}{%
  \textup{2010} Mathematics Subject Classification}
\makeatother
 \vspace{4mm}

\begin{document}
 \baselineskip=17pt
\hbox{Acta Arith. 183(2018), no.\,4, 339--356.}
\medskip

\title[Some variants of Lagrange's four squares theorem]{Some variants of Lagrange's four squares theorem}

\author[Y.-C. Sun]{Yu-Chen Sun}
\address{(Yu-Chen Sun) Department of Mathematics\\ Nanjing University\\
Nanjing 210093, People's Republic of China}
\email{syc@smail.nju.edu.cn}

\author[Z.-W. Sun]{Zhi-Wei Sun}
\address{(Zhi-Wei Sun, corresponding author) Department of Mathematics\\ Nanjing University\\
Nanjing 210093, People's Republic of China}
\email{zwsun@nju.edu.cn}

\date{}

\begin{abstract}
Lagrange's four squares theorem is a classical theorem in number theory.
Recently, Z.-W. Sun found that it can be further refined in various ways. In this paper we study some conjectures of Sun and obtain various refinements
of Lagrange's theorem.
We show that any nonnegative integer can be written as $x^2+y^2+z^2+w^2$ $(x,y,z,w\in\Z)$ with $x+y+z+w$ (or $x+y+z+2w$, or $x+2y+3z+w$) a square (or a cube).
Also, every $n=0,1,2,\ldots$ can be represented by $x^2+y^2+z^2+w^2$ $(x,y,z,w\in\Z)$ with $x+y+3z$ (or $x+2y+3z$) a square (or a cube), and
each $n=0,1,2,\ldots$ can be written as $x^2+y^2+z^2+w^2$ $(x,y,z,w\in\Z)$ with $(10w+5x)^2+(12y+36z)^2$ (or $x^2y^2+9y^2z^2+9z^2x^2$) a square.
We also provide an advance on the 1-3-5 conjecture of Sun. Our main results are proved by a new approach involving Euler's four-square identity.
\end{abstract}



\maketitle

\section{Introduction}
\setcounter{lemma}{0}
\setcounter{theorem}{0}
\setcounter{corollary}{0}
\setcounter{remark}{0}
\setcounter{equation}{0}
\setcounter{conjecture}{0}

In 1770 Lagrange showed that any $n\in\N=\{0,1,2,\ldots\}$ can be written as the sum of four squares.
This celebrated result is now known as Lagrange's four squares theorem.
The classical proof of this theorem (cf. \cite[pp.\,5-7]{N} or \cite{MW})
depends heavily on Euler's four-square identity. It is known that (cf. \cite{G}) only the following numbers have a unique representation as the sum of four unordered squares:
$$1,\ 3,\ 5,\ 7,\ 11,\ 15,\ 23,\ 2^{2k+1}m\ (k\in\N\ \t{and}\ m = 1,3,7).$$

In 1834, Jacobi (cf. \cite[p.\,59]{B} or \cite[p.\,83]{W}) considered the fourth power of the theta function
$$\varphi(q)=\sum_{n=-\infty}^\infty q^{n^2}\ \ (|q|<1)$$
and this led him to show that
$$r_4(n)=8\sum_{d\mid n\ \&\ 4\nmid d}d\quad\t{for all}\ n\in\Z^+=\{1,2,3,\ldots\},$$
where
$$r_4(n):=|\{(w,x,y,z)\in\Z^4:\ w^2+x^2+y^2+z^2=n\}|.$$

Recently, Sun \cite{S} found that Lagrange's four squares theorem can be refined in various surprising ways.
He \cite{S} called a polynomial $P(x,y,z,w)$ with integer coefficients {\it suitable} if any $n\in\N$ can be written as
$x^2+y^2+z^2+w^2$ with $x,y,z,w\in\N$ such that $P(x,y,z,w)$ a square. He proved that
$$x,\ 2x,\ x-y,\ 2x-2y,\ x^2y^2+y^2z^2+z^2x^2,\ x^2y^2+4y^2z^2+4z^2x^2$$
are all suitable, and conjectured that there are only finitely many linear suitable polynomials $ax+by+cz+dw$
with $a,b,c,d$ integers and $\gcd(a,b,c,d)$ squarefree (cf. the comments following \cite[Definition 1.1]{S}, and \cite[Conjectures 4.1, 4.3-4.4 and 4.12]{S}). Using the theory of ternary quadratic forms, he \cite{S} showed that
any $n\in\N$ can be written as $x^2+y^2+z^2+w^2$ with $x,y,z,w\in\Z$ such that $x+y+z$ (or $x+y+2z$) is a square (or a cube).

Our first theorem was originally conjectured by Sun \cite[Conjecture 4.1 and Remark 4.10]{S}.

\begin{theorem}\label{Th1.1} {\rm (i)} The polynomial $x^2+8yz$ is suitable, i.e., any $n\in\N$ can be written as $x^2+y^2+z^2+w^2$
$(x,y,z,w\in\N)$ with $x^2+8yz$ a square.

{\rm (ii)} The polynomial $x+2y$ is suitable. Moreover, any $n\in\Z^+$ can be written as $x^2+y^2+z^2+w^2$ $(x,y,z,w\in\N)$ with $x+2y$ a positive square.
\end{theorem}
\begin{remark}\label{Rem1.1} In contrast with Theorem 1.1(ii), Sun \cite[Theorem 1.1(iii)]{S} proved that any $n\in\N$ can be written as $x^2+y^2+z^2+w^2\ (x,y,z,w\in\Z)$
with $x+2y$ a square (or a cube).  Sun \cite[Conjecture 4.1]{S} also conjectured that $2x-y$ is suitable, but we are unable to show this.
\end{remark}

Our following results deal with integer variables.

\begin{theorem}\label{Th1.2} {\rm (i)} If $P(x,y,z)$ is one of the polynomials $x-y$, $x-2y$, $x+y-z$ and $x+y-2z$, then
any $n\in\N$ can be written as $x^2+y^2+z^2+w^2$ with $x,y,z,w\in\Z$ and $P(x,y,z)(x+2y+2z)=0$.

{\rm (ii)} Let $d\in\{1,2,3\}$ and $m\in\{2,3\}$. Then any $n\in\N$ can be written as $x^2+y^2+z^2+w^2$ with $x,y,z,w\in\Z$ such that $x+2y+2z=dt^m$ for some $t\in\Z$.
\end{theorem}
\begin{remark}\label{Rem1.2} Sun \cite[Remark 1.4, parts (iii) and (iv) of Conjecture 4.3]{S} conjectured that the polynomials
\begin{gather*}(x-y)(x+2y-2z),\ (x-2y)(x+2y-2z),
\\ (x+y-z)(x+2y-2z),\ 2x-y-2z,\ 2x+y-2z\end{gather*}
are all suitable.
\end{remark}

\begin{theorem}\label{Th1.3} {\rm (i)} For each $m=2,3$, any $n\in\N$ can be written as $x^2+y^2+z^2+w^2$ with $x,y,z,w\in\Z$ such that $x+y+z+w$ is an $m$th power.

{\rm (ii)} Let $m,n\in\Z^+$ with $n\eq 28\pmod{32}$. If $n$ can be written as $x^2+y^2+z^2+w^2$ with $x,y,z,w\in\Z$ such that $x+y+z+w$ is an $m$-th power,
then we must have $m\ls 3$.
\end{theorem}
\begin{remark}\label{Rem1.3} Sun \cite[Conjecture 4.4(v)]{S} conjectured that any positive integer can be written as $w^2+x^2+y^2+z^2$ with $w+x+y-z$ a square, where $w\in\Z$ and $x,y,z\in\N$
with $|w|\ls x\gs y\ls z<x+y$, and that any $n\in\N$ can be written as  $w^2+x^2 + y^2 + z^2$ with $w+x+y-z$ a nonnegative cube, where $w,x,y,z$ are integers with $|x|\ls y\gs z\gs0$.
\end{remark}

{\it Example} 1.1. We illustrate Theorem \ref{Th1.3}(i) with $m=2$ by an example:
$$60=1^2+3^2+5^2+(-5)^2\quad\t{with}\ \ 1+3+5+(-5)=2^2.$$

 We establish the following general theorem via a new approach involving Euler's four-square identity.

\begin{theorem}\label{Th1.4} Let $a,b,c,d\in\Z$ with $a,b,c,d$ not all zero, and let $\lambda\in\{1,2\}$ and $m\in\{2,3\}$.
Then any $n\in\N$ can be written as $x^2+y^2+z^2+w^2$ with $x,y,z,w\in\Z/(a^2+b^2+c^2+d^2)$ such that
$ax+by+cz+dw=\lambda r^m$ for some $r\in\N$.
\end{theorem}

With the help of Theorem \ref{Th1.4}, we obtain the following results.

\begin{theorem}\label{Th1.5} {\rm (i)} Let $\lambda\in\{1,2\}$ and $m\in\{2,3\}$.
Then any $n\in\N$ can be written as $x^2+y^2+z^2+w^2\ (x,y,z,w\in\Z)$ with $x+2y+3z-w\in\{\lambda r^m:\ r\in\N\}$.

{\rm (ii)} Any $n\in\N$ can be written as $w^2+x^2+y^2+z^2$ with $w,x,y,z\in\Z$
for which $w(x+2y+3z-w)=0$ and hence $w(x+2y+3z)$ is a square.

{\rm (iii)} For any $n\in\N$, there are $w,x,y,z\in\Z$ with $w^2+x^2+y^2+z^2=n$ for which
$(2w+x)(2w+x-y-3z)=0$ and hence $(10w+5x)^2+(12y+36z)^2$ is a square.
\end{theorem}
\begin{remark}\label{Rem1.4}. In contrast with Theorem \ref{Th1.5}(i), Sun \cite[Conjectures 4.4(iv) and 4.12(i)]{S} conjectured that $x+2y+3z-w$ and $x-2y+3z-w$ are both suitable.
Concerning parts (ii) and (iii) of Theorem \ref{Th1.5}, Sun \cite[Conjectures 4.5 and 4.8(i)]{S} also conjectured that both $w(x+2y+3z)$ and $(10w+5x)^2+(12y+36z)^2$ are suitable.
\end{remark}

\begin{theorem}\label{Th1.6} {\rm (i)} Any $n\in\N$ can be written as $w^2+x^2+y^2+z^2\ (w,x,y,z\in\Z)$ with $w(x+y+3z)=0$ $($or $(x-y)(x+y+3z)=0)$.
Consequently, each $n\in\N$ can be written as $w^2+x^2+y^2+z^2\ (w,x,y,z\in\Z)$ with $(x+y)^2+(4z)^2$ a square.

{\rm (ii)} Let $\lambda\in\{1,2\}$ and $m\in\{2,3\}$. Then any $n\in\N$ can be written as $x^2+y^2+z^2+w^2\ (x,y,z,w\in\Z)$ with $x+y+3z\in\{\lambda r^m:\ r\in\N\}$.
\end{theorem}
\begin{remark}\label{Rem1.5} Sun \cite[Conjecture 4.3(iv), Remark 1.4 and Conjecture 4.8(ii)]{S} conjectured that $x-y+3z$ and $(x-y)(x+y-3z)$ are suitable, and that
any positive integer can be written as $x^2+y^2+z^2+w^2$ with $x,y,z,w\in\N$, $y>z$, and $(x+y)^2+(4z)^2$ a square. We are also able to
show Theorem 1.6(iii) with $x,z\in\N$.
\end{remark}

\begin{corollary}\label{Cor1.1} Let $P(x,y,z)$ be one of the polynomials
$$x^2y^2+9y^2z^2+9z^2x^2,\ x^2+12yz,\ 9x^2-4yz,\ (x+3y)z.$$
Then any $n\in\Z^+$ can be written as $w^2+x^2+y^2+z^2$
$(w\in\Z\sm\{0\}$ and $x,y,z\in\Z)$ with $P(x,y,z)$ a square.
\end{corollary}
\begin{remark}\label{Rem1.6} Sun \cite[Remark 1.4 and Conjecture 4.10(i)]{S} conjectured that the polynomials $x^2y^2+9y^2z^2+9z^2x^2$, $x^2+kyz\ (k=12,24,32,48,84,120,252)$ and $9x^2-4yz$ are all suitable.
He also conjectured that any $n\in\Z^+$ can be written as $w^2+x^2+y^2+z^2$ with $w,x,y\in\N$ and $z\in\Z^+$ such that $(x+3y)z$ is a square (cf. \cite[Conjecture 4.2(i)]{S}).
\end{remark}

\begin{theorem}\label{Th1.7} {\rm (i)} Each $n\in\N$ can be written as $x^2+y^2+z^2+w^2\ (x,y,z,w\in\Z)$ with $w(x+y+z+2w)=0$.

{\rm (ii)} Let $\lambda\in\{1,2\}$ and $m\in\{2,3\}$. Then any $n\in\N$ can be written as $x^2+y^2+z^2+w^2\ (x,y,z,w\in\Z)$ with $x+y+z+2w=\lambda r^m$ for some $r\in\N$.

{\rm (iii)} Any $n\in\N$ can be written as $x^2+y^2+z^2+w^2\ (x,y,z,w\in\Z)$ with $w(x+2y+3z)=0$.

{\rm (iv)} For each $\lambda=1,2,4$ and $m\in\{2,3\}$, any $n\in\N$ can be written as $x^2+y^2+z^2+w^2\ (x,y,z,w\in\Z)$ with $x+2y+3z\in\{\lambda r^m:\ r\in\N\}$.
\end{theorem}
\begin{remark}\ \label{Rem1.7} Sun \cite[Conjectures 4.4(iv) and 1.2(ii)]{S} conjectured that both $x+y-z+2w$ and $x-y-z+2w$ are suitable;
in addition to Theorem \ref{Th1.7}(ii), we are also able to show that any $n\in\N$ can be written as $x^2+y^2+z^2+w^2$ with $x,y,w\in\N$ and $z\in\Z$
such that $x+y+z+2w$ is a square. Theorem \ref{Th1.7}(iv) provides an advance on Sun's conjecture that both $x+2y-3z$ and $x-2y+3z$ are suitable
(cf. \cite[parts (iii) and (iv) of Conjecture 4.3]{S}).
\end{remark}

\begin{corollary}\label{Cor1.2} Any $n\in\Z^+$ can be written as $x^2+y^2+z^2+w^2$
$(x,y,z\in\Z$ and $w\in\Z\sm\{0\})$ with $x^2+24yz$ a square. Also, each $n\in\Z^+$ can be written as $x^2+y^2+z^2+w^2$
$(x,y,z\in\Z$ and $w\in\Z\sm\{0\})$ with $4x^2y^2+9y^2z^2+36z^2x^2$ a square.
\end{corollary}

\begin{theorem}\label{Th1.8} Let $\lambda\in\{1,2\}$, $m\in\{2,3\}$ and $n\in\N$.

{\rm (i)} We can write $n$ as $x^2+y^2+z^2+w^2$ with $x,y,5z,5w\in\Z$ such that $x+3y+5z\in\{\lambda r^m:\ r\in\N\}$.

{\rm (ii)} We can write $n$ as $x^2+y^2+z^2+w^2$ with $x,y,z,w\in\Z/7$ such that $x+3y+5z\in\{\lambda r^m:\ r\in\N\}$.
\end{theorem}
\begin{remark}\label{Rem1.8} Theorem \ref{Th1.8} provides a remarkable progress on the 1-3-5 conjecture of Sun (cf. \cite[Conjecture 4.3(i)]{S}) which asserts that $x+3y+5z$ is suitable.
Sun \cite[Remark 4.3]{S} also conjectured any $n\in\N$ can be written as $x^2+y^2+z^2+w^2\ (x,y,z,w\in\Z)$ with $x+3y+5z$ a cube.
\end{remark}

We are going to show Theorems 1.1-1.3, 1.4-1.6 and 1.7-1.8 in Sections 2, 3 and 4 respectively.

Throughout this paper, for a nonzero integer $m$ and a prime power $p^a$ with $p$ prime and $a\in\N$, by $p^a\|m$ we mean $p^a\mid m$ but $p^{a+1}\nmid m$.

\maketitle

\section{Proofs of Theorems 1.1-1.3}
\setcounter{lemma}{0}
\setcounter{theorem}{0}
\setcounter{corollary}{0}
\setcounter{remark}{0}
\setcounter{equation}{0}
\setcounter{conjecture}{0}

For any $a,b,c\in\Z^+=\{1,2,3,\ldots\}$, we define
$$E(a,b,c):=\N\sm\{ax^2+by^2+cz^2: x,y,z\in\Z\}.$$

We need the following well-known theorem.
\medskip

\noindent{\bf Gauss-Legendre Theorem} {\rm (\cite[p.\,23]{N})}. We have
\begin{equation}\label{2.1}E(1,1,1)=\{4^k(8l+7):\ k,l\in\N\}.\end{equation}

\medskip
\noindent {\it Proof of Theorem} \ref{Th1.1}. (i) If $n\not\in E(1,1,1)$, then for some $w,x,y\in\N$
we have $n=w^2+x^2+y^2+0^2$ with $x^2+8y\times0=x^2$.

Now suppose that $n\in E(1,1,1)$. As $E(1,2,3)\cap E(1,1,1)=\emptyset$ by \cite[Lemma 3.2]{S}, we can write $n$ as $w^2+2v^2+3z^2$
with $w,v,z\in\N$. Let $x=|v-z|$ and $y=v+z$. Then
$$n=w^2+(v-z)^2+(v+z)^2+z^2=w^2+x^2+y^2+z^2$$
and
$$x^2+8yz=|v-z|^2+8(v+z)z=(v+3z)^2.$$
This proves Theorem 1.1(i).

(ii) Note that
$$\l\lfloor\l(\f{10}{\root4\of{5}-1}\r)^4\r\rfloor=166094.$$
For $n=1,2,\ldots, 166094$ we can verify the desired result via a computer.

Now fix $n\in\N$ with $n>166094$, and assume that any $m=0,\ldots,n-1$ can be written as $x^2+y^2+z^2+w^2$ $(x,y,z,w\in\N$)
with $x+2y$ a positive square.

If $16\mid n$, then by the induction hypothesis there are $x,y,z,w\in\N$ such that
$x^2+y^2+z^2+w^2=n/16$ and $x+2y$ is a positive square, hence
$n=(4x)^2+(4y)^2+(4z)^2+(4w)^2$ with $4x+2(4y)=2^2(x+2y)$ a positive square.

Below we suppose that $16\nmid n$. As
$n>166094$, we have $\root4\of{5n}-\root4\of{n}\gs 10$
and hence for each $r=0,1$ there is a positive integer $j_r$ such that
$$\root4\of{n}\ls 5(2j_r+r)<\root4\of{5n}.$$
We claim that $5n-(5(2j_r+r))^4\not\in E(1,1,1)$ for some $r=0,1$.

{\it Case} 1. $2\nmid n$.

If $5n-(10j_0)^4\eq7\pmod 8$, then $5n-(10j_1+5)^4\eq6\pmod 8$.
So, $5n-(10j_0)^4$ or $5n-(10j_1+5)^4$ does not belong to $E(1,1,1)$.

{\it Case} 2. $4\|n$.

Write $5n=4q$ with $q$ odd. If $q\eq1\pmod 4$, then
$5n-(10j_0)^4\eq4\pmod{16}$ and hence $5n-(10j_0)^4\not\in E(1,1,1)$.
If $q\eq3\pmod 4$, then $5n-(10j_1+5)^4\eq3\pmod 8$ and hence $5n-(10j_1+5)^4\not\in E(1,1,1)$.

{\it Case} 3. $2\|n$ or $8\|n$.

In this case, $2\|5n-(10j_0)^4$ or $8\|5n-(10j_0)^4$, and hence $5n-(10j_0)^4\not\in E(1,1,1)$.

By the claim, for some $j\in\{2j_0,2j_1+1\}$ we can write $5n-(5j)^4=t^2+u^2+v^2$ with $t,u,v\in\N$.
As a square is congruent to one of $0,1,-1$ modulo $5$, one of $t,u,v$ must be divisible by $5$.
Without loss of generality, we assume that $5\mid t$ and $u\eq2v\pmod 5$.
(Note that $5\mid u^2+v^2\iff u^2\eq(2v)^2\pmod 5$.)
Set $s=(5j)^2$. Then
\begin{align*} n=&\f{s^2+t^2}5+\f{u^2+v^2}5=5\l(\f s5\r)^2+5\l(\f t5\r)^2+\f{u^2+v^2}5
\\=&\l(\f{s+2t}5\r)^2+\l(\f{2s-t}5\r)^2+\l(\f{u-2v}5\r)^2+\l(\f{2u+v}5\r)^2
\end{align*}
with
$$\f{s+2t}5+2\f{2s-t}5=s=(5j)^2.$$
As $s^2\gs n$, we have $t^2\ls 5n-s^2\ls4s^2$
and hence $2s-t\gs0$. This concludes our induction proof. \qed

\begin{lemma}\label{Lem2.1} Suppose that $n$ is the sum of three squares. Then we can write
$n=x^2+y^2+z^2+w^2$ $(x,y,z,w\in\Z)$ with $x+2y+2z=0$.
\end{lemma}
\Proof. Clearly, $n=u^2+v^2+w^2$ for some integers $u,v,w\not\eq2\pmod 3$.
(Note that if $x\eq2\pmod 3$ then $-x\not\eq2\pmod 3$.) Obviously, two of $u,v,w$ are congruent modulo $3$.
Without loss of generality, we suppose that $u\eq v\pmod 3$. Set
$$y=\f{2u+v}3\ \ \t{and}\ \ z=\f{u+2v}3.$$
Then
$$n=u^2+v^2+w^2=(2y-z)^2+(2z-y)^2+w^2=(2y-2z)^2+(-y)^2+z^2+w^2$$
with $(2y-2z)+2(-y)+2z=0$.
This ends the proof. \qed

\medskip
\noindent {\it Proof of Theorem} \ref{Th1.2}. (i) By \cite[Lemma 3.2]{S}, all
the sets
$$E(1,1,2),\ E(1,1,5),\ E(1,2,3),\ E(1,2,6)$$
are subsets of $\N\sm E(1,1,1)$. If $n\not\in E(1,1,2)$, then we can write
$n=x^2+y^2+z^2+w^2$ with $x,y,z,w\in\Z$ and $x-y=0$. If $n\not\in E(1,1,5)$, then
we can write $n=x^2+y^2+z^2+w^2$ with $x,y,z,w\in\Z$ and $x-2y=0$.
When $n\not\in E(1,2,6)$, by \cite[Lemma 3.3]{S} we can write $n=x^2+y^2+z^2+w^2$ with $x,y,z,w\in\N$ and $x+y-z=0$.
If $n\not\in E(1,2,3)$, then again by \cite[Lemma 3.3]{S} we can write $n=x^2+y^2+z^2+w^2$ with $x,y,z,w\in\Z$ and $x+y-2z=0$.
Combining these with Lemma 2.1 we immediately finish the proof of Theorem 1.2(i).

(ii) For $n=0,1,\ldots,4^m-1$, we can verify the desired result in Theorem 1.2(ii) directly.

Now fix an integer $n\gs 4^m$ and assume the desired result for smaller values of $n$.

{\it Case} 1. $4^m\mid n$.

In this case, by the induction hypothesis, there are $x,y,z,w\in\Z$ with $n/4^m=x^2+y^2+z^2+w^2$ such that $x+2y+2z=dt^m$ for some $t\in\Z$.
Hence $n=(2^mx)^2+(2^my)^2+(2^mz)^2+(2^mw)^2$ with $2^mx+2(2^my)+2(2^mz)=2^m(x+2y+2z)=d(2t)^m$.
\medskip

{\it Case} 2. $4^m\nmid n$ and $n\not=4^k(8l+7)$ for any $k,l\in\N$.

In this case, by the Gauss-Legendre theorem and Lemma \ref{Lem2.1}, we can write $n$ as $x^2+y^2+z^2+w^2$
$(x,y,z,w\in\Z$) with $x+2y+2z=0^m$.
\medskip

{\it Case} 3. $n=4^k(8l+7)$ for some $k,l\in\N$ with $k<m$.

We first consider the subcase $d\in\{1,2\}$. We claim that $\{9n-d^2,9n-d^24^m\}\not\se \{4^r(8q+7):\ q,r\in\N\}$.
In fact, for $d=2$ we have
$$9n-2^2=9\times4^k(8l+7)-4=\begin{cases}8(9l+7)+3&\t{if}\ k=0;
\\4(8(9l+7)+6)&\t{if}\ k=1,
\\4(8(36l+31)+3)&\t{if}\ k=2.
\end{cases}$$
For $d=1$, clearly
$$9n-1^2=9\times 4^k(8l+7)-1\eq \begin{cases} 6\pmod 8&\t{if}\ k=0,
\\3\pmod 8&\t{if}\ k=1,
\end{cases}$$
and when $k=2<m=3$ we have
$$9n-4^3=9\times 4^2(8l+7)-4^3=4^2(72l+59)=4^2(8(9l+7)+3).$$
Thus, by the Gauss-Legendre theorem, for some $\da\in\{1,2^m\}$ there exist $a,b,c\in\Z$ such that $9n-(d\da)^2=a^2+b^2+c^2$. Clearly, we cannot have $3\nmid abc$ since $9n-(d\da)^2\not\eq0\pmod 3$.
Without loss of generality, we assume that $c=3w$ with $w\in\Z$. As $a^2+b^2\eq-(d\da)^2\eq2\pmod 3$, we must have $3\nmid ab$.
We may simply suppose that $a=3u+2d\da$ and $b=3v-2d\da$ with $u,v\in\Z$. (Note that if $x\eq1\eq-2\pmod3$ then $-x\eq2\pmod 3$.)
Since
\begin{align*} &12d\da u-12d\da v+8(d\da)^2
\\\eq&(3u+2d\da)^2+(3v-2d\da)^2=a^2+b^2
\eq-(d\da)^2\pmod 9,
\end{align*}
we must have $3\mid d\da(u-v)$ and hence $u\eq v\pmod 3$.
Set
$$y=\f{2u+v}3\ \ \t{and}\ \ z=\f{u+2v}3.$$
Then
\begin{align*} 9n-(d\da)^2=&a^2+b^2+c^2=(3u+2d\da)^2+(3v-2d\da)^2+9w^2
\\=&(3(2y-z)+2d\da)^2+(3(2z-y)-2d\da)^2+9w^2
\\=&9(2y-z)^2+9(2z-y)^2+12d\da((2y-z)-(2z-y))
\\&+8(d\da)^2+9w^2
\end{align*}
and hence
\begin{align*}n=&(2y-z)^2+(2z-y)^2+4d\da(y-z)+w^2+(d\da)^2
\\=&(2y-2z+d\da)^2+(-y)^2+z^2+w^2
\end{align*}
with $(2y-2z+d\da)+2(-y)+2z=d\da\in\{dt^m:\ t=1,2\}$.

Now we consider the subcase $d=3$. For $q=0,\ldots,363$, we can verify via a computer that $4^2(8q+7)$
can be written as $x^2+y^2+z^2+w^2$ with $x,y,z,w\in\Z$ such that $x+2y+2z=3t^3$
for some $t\in\Z$; for example,
$$4^2(8\times362+7)=46448=8^2+(-148)^2+156^2+12^2$$
with $8+2(-148)+2\times156=3\times2^3$. From now on we assume that $k<2$ or $l\gs 364$.
When $k\in\{0,1\}$, it is easy to see that
$n-3^{2m}$ is congruent to $6$ or $3$ modulo $8$. If $m=3$ and $k=2$, then
$$n-6^{2m}=4^2(8l+7)-4^3\times 3^6=4^2(8(l-364)+3).$$
So, for some $\da\in\{3^m,6^m\}$ we have $n-\da^2\not\in\{4^i(8j+7):\ i,j\in\N\}$.
By the Gauss-Legendre theorem, $n-\da^2=a^2+b^2+c^2$ for some $a,b,c\in\Z$.
Clearly, two of $a^2,b^2,c^2$, say $a^2$ and $b^2$, are congruent modulo $3$.
Without loss of generality, we may assume that $a\eq b\pmod 3$. Let $u=a-2\da$ and $v=b+2\da$.
Then $u\eq v\pmod 3$ since $3\mid \da$.
Set
$$y=\f{2u+v}3\ \ \t{and}\ \ z=\f{u+2v}3.$$
Then
$$n-\da^2=(u+2\da)^2+(v-2\da)^2+c^2=(2y-z+2\da)^2+(2z-y-2\da)^2+c^2.$$
It follows that
$$n=(2y-2z+3\da)^2+(-y)^2+z^2+c^2$$
with $(2y-2z+3\da)+2(-y)+2z=3\da\in\{3\times t^m:\ t=3,6\}$ as desired.

\medskip

So far we have completed the proof of Theorem \ref{Th1.2}. \qed

\medskip
\noindent {\it Proof of Theorem} \ref{Th1.3}.
(i) For $n=0,1,\ldots,4^m-1$, the desired result can be verified directly.

Now let $n\gs 4^m$ be an integer and assume the desired result for smaller values of $n$.

{\it Case}\ 1. $4^m\mid n$.

By the induction hypothesis, we can write
$$\f n{4^m}=x^2+y^2+z^2+w^2\ \ \t{with}\ x,y,z,w\in\Z$$
such that $x+y+z+w=s^m$ for some $s\in\Z$. Therefore,
$$n=(2^mx)^2+(2^my)^2+(2^mz)^2+(2^mw)^2$$
with $2^mx+2^my+2^mz+2^mw=(2s)^m$.
\medskip

{\it Case} 2. $2\mid n$ but $4^m\nmid n$.

If $n\not=4^k(8l+7)$ for all $k,l\in\N$ with $k<m$, then by the Gauss-Legendre theorem we can write $n=u^2+v^2+w^2$ with $u,v,w\in\Z$.
Clearly, $u+v+w\eq n\eq0\pmod 2$.
Let
$$t=\f{u+v+w}2,\ x=\f{-u+v+w}2,\ y=\f{u-v+w}2\ \t{and}\ z =\f{u+v-w}2.$$
Then $x+y=w$, $y+z=u$, $x+z=v$, and $x+y+z=t$. Therefore,
\begin{align*} n=&u^2+v^2+w^2=(y+z)^2+(x+z)^2+(x+y)^2
\\=&(x+y+z)^2+x^2+y^2+z^2=(-t)^2+x^2+y^2+z^2
\end{align*}
with $(-t)+x+y+z=0^m$.

When $n=4^k(8l+7)$ for some $k,l\in\N$ with $k<m$, we have $k\in\{1,2\}$ and hence $n-4^{m-1}\not=4^a(8b+7)$ for any $a,b\in\N$.
Thus, by the Gauss-Legendre theorem
there are $u,v,w\in\Z$ such that $n-4^{m-1}=(u-2^{m-1})^2+(v-2^{m-1})^2+(w-2^{m-1})^2$. Set
$$t=\f{u+v+w}2,\ x=\f{-u+v+w}2,\ y=\f{u-v+w}2\ \t{and}\ z =\f{u+v-w}2.$$
Then
\begin{align*} n=&(u-2^{m-1})^2+(v-2^{m-1})^2+(w-2^{m-1})^2+4^{m-1}
\\=&(y+z-2^{m-1})^2+(x+z-2^{m-1})^2+(x+y-2^{m-1})^2+4^{m-1}
\\=&(x+y+z-2^m)^2+x^2+y^2+z^2=(2^m-t)^2+x^2+y^2+z^2
\end{align*}
with $(2^m-t)+x+y+z=2^m$.
\medskip

{\it Case} 3. $2\nmid n$.

As $4n-1\eq3\pmod 8$, by the Gauss-Legendre theorem there are $u,v,w\in\Z$ such that
$$4n-1=(2u-1)^2+(2v-1)^2+(2w-1)^2=(2u'-1)^2+(2v-1)^2+(2w-1)^2,$$
where $u'=1-u$.

When $u+v+w$ is even, we set
$$t=\f{u+v+w}2,\ x=\f{-u+v+w}2,\ y=\f{u-v+w}2\ \t{and}\ z =\f{u+v-w}2,$$
and deduce that
\begin{align*} 4n-1=&(2u-1)^2+(2v-1)^2+(2w-1)^2
\\=&(2(y+z)-1)^2+(2(x+z)-1)^2+(2(x+y)-1)^2
\\=&4\l((x+y+z-1)^2+x^2+y^2+z^2\r)-1.
\end{align*}
Therefore
$$n=(1-t)^2+x^2+y^2+z^2\ \t{with}\ (1-t)+x+y+z=1^m.$$

Now assume that $u+v+w$ is odd. Then $u'+v+w$ is even. Set
$$t'=\f{u'+v+w}2,\ x=\f{-u'+v+w}2,\ y=\f{u'-v+w}2\ \t{and}\ z =\f{u'+v-w}2,$$
Then we get
\begin{align*} 4n-1=&(2u'-1)^2+(2v-1)^2+(2w-1)^2
\\=&(2(y+z)-1)^2+(2(x+z)-1)^2+(2(x+y)-1)^2
\end{align*}
and hence
$$n=(1-t')^2+x^2+y^2+z^2\ \ \ \t{with}\ (1-t')+x+y+z=1^m.$$
This concludes the proof of the first part of Theorem 1.3.

(ii) Write $n=4(8l+7)$ with $l\in\N$. Suppose that
$n$ can be written as $x^2+y^2+z^2+w^2$ with $x,y,z,w\in\Z$
such that $x+y+z+w=t^m$ for some $t\in\Z$. Then
\begin{align*} n=&(x+y+z-t^m)^2+x^2+y^2+z^2
\\=&(x+y)^2+(y+z)^2+(z+x)^2-t^m((x+y)+(y+z)+(z+x))+t^{2m}
\end{align*}
and hence
$$4n-t^{2m}=(2(x+y)-t^m)^2+(2(y+z)-t^m)^2+(2(z+x)-t^m)^2$$
If $m\gs 4$ and $t=2s$ with $s\in\Z$, then
$$4n-t^{2m}=4^2(8l+7)-(2s)^{2m}=4^2(8(l-2^{2m-7}s^{2m})+7).$$
If $m\gs 4$ and $2\nmid t$, then
$$4n-t^{2m}=4^2(8l+7)-t^{2m}\eq -1\eq7\pmod 8.$$
So, in light of the Gauss-Legendre theorem, we must have $m\ls 3$.

The proof of Theorem \ref{Th1.3} is now complete. \qed

\maketitle

\section{Proofs of Theorems 1.4-1.6 and Corollary 1.1}
\setcounter{lemma}{0}
\setcounter{theorem}{0}
\setcounter{corollary}{0}
\setcounter{remark}{0}
\setcounter{equation}{0}
\setcounter{conjecture}{0}

\begin{lemma}\label{Lem3.1} Let $a,b,c,d$ be integers not all zero, and let
\begin{equation}\label{3.1}\begin{cases} x=\f{as-bt-cu-dv}{a^2+b^2+c^2+d^2},
\\y=\f{bs+at+du-cv}{a^2+b^2+c^2+d^2},
\\z=\f{cs-dt+au+bv}{a^2+b^2+c^2+d^2},
\\w=\f{ds+ct-bu+av}{a^2+b^2+c^2+d^2}.
\end{cases}\end{equation}
Then
\begin{equation}\label{3.2}\begin{cases} ax+by+cz+dw=s,
\\ay-bx+cw-dz=t,
\\az-bw-cx+dy=u,
\\aw+bz-cy-dx=v,
\end{cases}
\end{equation}
and
\begin{equation}\label{3.3}(a^2+b^2+c^2+d^2)(x^2+y^2+z^2+w^2)=s^2+t^2+u^2+v^2.\end{equation}
\end{lemma}
\Proof. It is easy to verify (\ref{3.2}) directly. By (\ref{3.2}) and Euler's four-square identity, we immediately have (\ref{3.3}).  \qed

\medskip
\noindent{\it Proof of Theorem} \ref{Th1.4}. Let $n\in\N$. By \cite[Theorem 1.1]{S},
we can write $(a^2+b^2+c^2+d^2)n$ as $(\lambda r^m)^2+t^2+u^2+v^2$ with $r,t,u,v\in\N$. Set $s=\lambda r^m$, and define $x,y,z,w$ by (\ref{3.1}).
Then $x,y,z,w\in\Z/(a^2+b^2+c^2+d^2)$. By Lemma \ref{Lem3.1}, we have
$$x^2+y^2+z^2+w^2=\f{s^2+t^2+u^2+v^2}{a^2+b^2+c^2+d^2}=n$$
and
$$ax+by+cz+dw=s=\lambda r^m.$$ \qed

\begin{lemma}\label{Lem3.2} Suppose that $3n=s^2+t^2+u^2+v^2$ with $n\in\N$ and $s,t,u,v\in\Z$. Define
\begin{equation}\label{3.4}\begin{cases} x=\f{s-t-u}3,\\y=\f{s+t-v}3,\\z=\f{s+u+v}3,\\w=\f{t-u+v}3.\end{cases}\end{equation}
Then we have
$$x^2+y^2+z^2+w^2=n\ \mbox{and}\ \ x+2y+3z-w=2s+u.$$
\end{lemma}
\Proof. Taking $a=b=c=1$ and $d=0$ in Lemma \ref{Lem3.1} we obtain that
$$x^2+y^2+z^2+w^2=\f{s^2+t^2+u^2+v^2}{1^2+1^2+1^2+0^2}=n.$$
The equality $x+2y+3z-w=2s+u$ follows from (\ref{3.4}).
\qed

\begin{lemma}\label{Lem3.3}  {\rm (Sun \cite[Theorem \ref{Th1.2}(iii)]{S})} Let $\lambda\in\{1,2\}$ and $m\in\{2,3\}$.
Then, any $n\in\N$ can be written as $x^2+y^2+z^2+w^2$ with $x,y,z,w\in\Z$ such that $x+2y=\lambda r^m$ for some $r\in\N$.
\end{lemma}

\medskip
\noindent{\it Proof of Theorem} \ref{Th1.5}. (i) By Lemma \ref{Lem3.3}, we can write $3n$ in the form $s^2+t^2+u^2+v^2$
with $s,t,u,v\in\Z$ and $2s+u=\lambda r^m$ for some $r\in\N$. Since $t^2+v^2\eq 2(s^2+u^2)=(s-u)^2+(s+u)^2\pmod3$,
and a square is congruent to $0$ or $1$ modulo 3, without loss of generality we may assume that
$t\eq s-u\pmod3$ and $v\eq-s-u\pmod 3$. Now define $x,y,z,w$ as in (\ref{3.4}). It is easy to see that $x,y,z,w\in\Z$.
In light of Lemma \ref{Lem3.2}, $x^2+y^2+z^2+w^2=n$ and $x+2y+3z-w=2s+u=\lambda r^m$. This proves Theorem \ref{Th1.5}(i).

(ii) If $n\not\in E(1,1,1)$, then for some $x,y,z\in\Z$ we have $n=0^2+x^2+y^2+z^2$ with $0(x+2y+3z-0)=0$ and $0(x+2y+3z)=0^2$.

Now suppose that $n\in E(1,1,1)$. By the Gauss-Legendre theorem, $n=4^k(8l+7)$ for some $k,l\in\N$.
Since $E(1,1,5)=\{4^i(8j+3):\ i,j\in\N\}$
by Dickson \cite[pp.\,112-113]{D}, we have $3n=4^k(8(3l+2)+5)\not\in E(1,1,5)$.
Hence there are $s,t,u,v\in\Z$ with $u=-2s$ such that $3n=5s^2+t^2+v^2=s^2+t^2+u^2+v^2$.
As $t^2+v^2\eq s^2\not\eq2\pmod3$, we have $3\mid tv$. Without loss of generality, we may assume that $3\mid t$ and $v\eq s\pmod 3$.
Let $x,y,z,w$ be given by (\ref{3.4}). Clearly, $x,y,z,w\in\Z$. In view of Lemma \ref{Lem3.2}, we have
$x^2+y^2+z^2+w^2=n$ and $x+2y+3z-w=2s+u=0$. Note that $w(x+2y+3z)=w^2$. This concludes our proof of Theorem \ref{Th1.5}(ii).

(iii) If $n\not\in E(1,1,5)$, then for some $w,x,y,z\in\Z$ we have $n=w^2+x^2+y^2+z^2$ with $x=-2w$, thus
$(2w+x)(2w+x-y-3z)=0$ and $(10w+5x)^2+(12y+36z)^2=(12y+36z)^2$.

Now suppose that $n\in E(1,1,5)$. By Dickson \cite[pp.\,112-113]{D}, $E(1,1,5)=\{4^k(8l+3):\ k,l\in\N\}$.
So, $n=4^k(8l+3)$ for some $k,l\in\N$. As $3n=4^k(8(3l+1)+1)\not\in E(1,1,5)$,
there are $s,t,u,v\in\Z$ with $u=-2s$ such that $3n=5s^2+t^2+v^2=s^2+t^2+u^2+v^2$.
As in the proof of Theorem \ref{Th1.5}(i), without loss of generality we may assume that $3\mid t$ and $v\eq s\pmod 3$.
Let $x,y,z,w$ be given by (\ref{3.4}). As in the proof of Theorem \ref{Th1.5}(ii), we have $x,y,z,w\in\Z$, $x^2+y^2+z^2+w^2=n$
and $x+2y+3z-w=2s+u=0$. Thus $n=y^2+x^2+(-z)^2+w^2$ with
$$(2y+x)(2y+x-w-3(-z))=0$$
and
$$(10y+5x)^2+(12w+36(-z))^2 = 5^2(2y+x)^2+12^2(w-3z)^2 = (13(2y+x))^2.$$
This proves Theorem \ref{Th1.5}(iii). \qed

\begin{lemma}\label{Lem3.4} Suppose that $11n=s^2+t^2+u^2+v^2$ with $n\in\N$ and $s,t,u,v\in\Z$. Define
\begin{equation}\label{3.5}\begin{cases} x=\f{s+t-3u}{11},\\y=\f{s-t+3v}{11},\\z=\f{3s+u-v}{11},\\w=\f{3t+u+v}{11}.\end{cases}
\end{equation}
Then we have
$$x^2+y^2+z^2+w^2=n\ \t{and}\ \ x+y+3z=s.$$
\end{lemma}
\Proof. Taking $a=1$, $b=-1$, $c=3$ and $d=0$ in Lemma \ref{Lem3.1} we obtain that
$$x^2+(-y)^2+z^2+w^2=\f{s^2+t^2+u^2+v^2}{1^2+(-1)^2+3^2+0^2}=n.$$
By (\ref{3.5}), we immediately see that $x+y+3z=s$. \qed

\begin{lemma}\label{Lem3.5} Let $m\in\{2,3\}$ and $n\in\N$.

{\rm (i)} Suppose that $4^m\nmid n$ and $n\gs\lambda^2(2q)^{2m}$, where $q$ is a positive odd integer.
 Then, for each $\lambda=1,2$, we can write $n$ as $(\lambda r^m)^2+x^2+y^2+z^2$ with $r\in\{0,q,2q\}$ and $x,y,z\in\Z$.

{\rm (ii)} Suppose that $4^m\nmid 2n$ and $2n\gs\lambda^2q^{2m}$, where $q$ is a positive odd integer.
 Then $2n$ can be written as $(4r^m)^2+x^2+y^2+z^2$ with $r\in\{0,q\}$ and $x,y,z\in\Z$.
\end{lemma}
\Proof.  (i) Suppose that $n$ is not the sum of three squares. By the Gauss-Legendre theorem, $n=4^k(8l+7)$ for some $k\in\{0,\ldots,m-1\}$ and $l\in\N$.
Note that $q^{2m}=8v+1$ for some $v\in\N$. Clearly,
$$n-\lambda^2q^{2m}=\begin{cases} 8l+7-\lambda^2(8v+1)\eq 3,6\pmod 8&\t{if}\ k=0,
\\4(8l+7)-\lambda^2(8v+1)\eq3\pmod 8&\t{if}\ k=\lambda=1,
\\4^k(8l-8v+6)&\t{if}\ k=1\ \text{and}\ \lambda=2,
\\4(8(4l-v+3)+3)&\t{if}\ k=\lambda=2.
\end{cases}$$
If $\lambda=1,\ k=2$ and $m=3$, then
$$n-\lambda^2(2q)^{2m}=4^2(8l+7)-2^6(8v+1)=4^2(8(l-4v)+3).$$
In view of the above and the Gauss-Legendre theorem, there exists $r\in\{0,q,2q\}$ such that $n-(\lambda r^m)^2=x^2+y^2+z^2$ for some $x,y,z\in\Z$.
This proves Lemma \ref{Lem3.5}(i).

(ii)  Suppose that $2n$ is not the sum of three squares. By the Gauss-Legendre theorem, $2n=4^k(8l+7)$ for some $k\in\{1,\ldots,m-1\}$ and $l\in\N$.
Note that $q^{2m}=8v+1$ for some $v\in\N$. Clearly,
$$2n-16q^{2m}=\begin{cases} 4(8l+7)-16(8v+1)=4(8(l+4v)+3)&\t{if}\ k=1,
\\4^2(8l+7)-16(8v+1)=4^2(8(l-v)+6)&\t{if}\ k=2.
\end{cases}$$

In light of the above and the Gauss-Legendre theorem, there exists $r\in\{0,q\}$ such that $2n-(4 r^m)^2=x^2+y^2+z^2$ for some $x,y,z\in\Z$.
This proves Lemma \ref{Lem3.5}(ii). \qed

\begin{lemma}\label{Lem3.6} {\rm (see, e.g., Dickson \cite[pp.\,112-113]{D})} We have
\begin{equation}\label{3.6}E(1,1,2)=\{4^k(16l+14):\ k,l\in\N\}.\end{equation}
\end{lemma}

\medskip
\noindent{\it Proof of Theorem} \ref{Th1.6}. (i) If $n\not\in E(1,1,1)$, then for some $x,y,z\in\Z$ we have $n=0^2+x^2+y^2+z^2$ with $0(x+y+3z)=0$,
 and also $(x+y)^2+(4\times0)^2=(x+y)^2$ is a square. If $n\not\in E(1,1,2)$, then for some $x,y,z,w\in\Z$ with $x=y$ we have $n=x^2+y^2+z^2+w^2$ with $(x-y)(x+y+3z)=0$,

If $n\in E(1,1,1)$, then by the Gauss-Legendre theorem, we have $n=4^k(8l+7)$ for some $k,l\in\N$, and hence $11n=4^k(8(11l+9)+5)\not\in E(1,1,1)$.
If $n\in E(1,1,2)$, then by (\ref{3.6}) $n=4^k(16l+14)$ for some $k,l\in\N$, and hence $11n=4^k(8(22l+19)+2)\not\in E(1,1,1)$ by the Gauss-Legendre theorem.

Now suppose that $11n\not\in E(1,1,1)$. We can write $11n=s^2+t^2+u^2+v^2$ with $s=0$ and $t,u,v\in\Z$.
Clearly, a square is congruent to one of $0,1,-2,3,5$ modulo $11$, and no sum of three distinct numbers among $0,1,-2,3,5$ is divisible by $11$.
Thus two of $t^2,u^2,v^2$ are congruent modulo $11$. Without loss of generality, we may assume that $u\eq v\pmod{11}$.
As $t^2\eq-u^2-v^2\eq(3u)^2\pmod{11}$, without loss of generality we simply assume that $t\eq3u\pmod{11}$.
Define $x,y,z,w$ as in (3.5). Then $x,y,z,w\in\Z$. By Lemma \ref{Lem3.4}, we have $x^2+y^2+z^2+w^2=n$ and $x+y+3z=s=0$.
Therefore $(x+y)^2+(4z)^2=(-3z)^2+(4z)^2=(5z)^2$.

(ii) For $n=0,1,\ldots,\lambda^2 22^{2m}$ we can verify the desired result via a computer.

Now fix an integer $n\gs \lambda^2 22^{2m}$ and assume that the desired result holds for smaller values of $n$.

If $4^m\mid n$, then by the induction hypothesis we can write $n/4^m=x^2+y^2+z^2+w^2\ (x,y,z,w\in\Z)$
with $x+y+3z=\lambda r^m$ for some $r\in\N$, hence
$$n=(2^mx)^2+(2^my)^2+(2^mz)^2+(2^mw)^2\ \ \t{with}\ 2^mx+2^my+3\times 2^mz=\lambda(2r)^m.$$

Below we suppose that $4^m\nmid n$. Applying Lemma \ref{Lem3.5}(i) with $q=11$, we can write $11n$ as $(\lambda r^m)^2+t^2+u^2+v^2$
with $r\in\{0,11,22\}$ and $t,u,v\in\Z$. As $s:=\lambda r^m\eq0\pmod{11}$, we have $t^2+u^2+v^2=11n-s^2\eq0\pmod{11}$.
Clearly, a square is congruent to one of $0,1,-2,3,4,5$ modulo $11$, and no sum of three distinct numbers among $0,1,-2,3,4,5$
is divisible by $11$. As $t^2+u^2+v^2\eq0\pmod{11}$, two of $t^2,u^2,v^2$ are congruent modulo $11$. Without loss of generality, we may assume that
$u\eq v\pmod{11}$ and $t\eq 3u\pmod{11}$. Define $x,y,z,w\in\Z$ by (\ref{3.5}). It is easy to see that $x,y,z,w\in\Z$.
By Lemma \ref{Lem3.4}, $x^2+y^2+z^2+w^2=n$ and $x+y+3z=s=\lambda r^m$. This ends our proof of Theorem \ref{Th1.6}(ii).

Combining the above, we have completed the proof of Theorem \ref{Th1.6}. \qed

\medskip
\noindent{\it Proof of Corollary} \ref{Cor1.1}. If $n\not\in E(1,1,1)$, then there are $w,x,y,z\in\Z$ with $w\not=0$ and $z=0$ such that
$n=w^2+x^2+y^2+z^2$, hence $x^2y^2+9y^2z^2+9z^2x^2=(xy)^2$, $(x+3y)z=0^2$, $x^2+12yz=x^2$ and $9x^2-4yz=(3x)^2$ are all squares.

Now suppose that $n\in E(1,1,1)$. By Theorem \ref{Th1.6}(i), $n$ can be written as $x^2+y^2+z^2+w^2$ with $w(x+y-3z)=0$.
Clearly, $w\not=0$ since $n\in E(1,1,1)$. Thus $x+y=3z$ and hence
\begin{align*} x^2y^2+9y^2z^2+9z^2x^2=&x^2y^2+(x^2+y^2)(3z)^2
\\=&x^2y^2+(x^2+xy+y^2-xy)(x^2+xy+y^2+xy)
\\=&(x^2+xy+y^2)^2.
\end{align*}
Note that
$$x^2+12yz=(3z-y)^2+12yz=(3z+y)^2$$
and $$9z^2-4xy=(x+y)^2-4xy=(x-y)^2.$$
Also, $n=w^2+(-x)^2+z^2+y^2$ with $(-x+3z)y=y^2$. This concludes our proof. \qed

\maketitle

\section{Proofs of Theorems 1.7-1.8 and Corollary 1.2}

\setcounter{lemma}{0}
\setcounter{theorem}{0}
\setcounter{corollary}{0}
\setcounter{remark}{0}
\setcounter{equation}{0}
\setcounter{conjecture}{0}

\begin{lemma}\label{Lem4.1} Suppose that $7n=s^2+t^2+u^2+v^2$ with $n\in\N$ and $s,t,u,v\in\Z$. Define
\begin{equation}\label{4.1}\begin{cases} x=\f{s-t-u-2v}7,\\y=\f{s+t+2u-v}7,\\ z=\f{s-2t+u+v}7,\\w=\f{2s+t-u+v}7.\end{cases}
\end{equation}
Then we have
$$x^2+y^2+z^2+w^2=n\ \t{and}\ \ x+y+z+2w=s.$$
\end{lemma}
\Proof. Taking $a=b=c=1$ and $d=2$ in Lemma \ref{Lem3.1} we obtain that
$$x^2+y^2+z^2+w^2=\f{s^2+t^2+u^2+v^2}{1^2+1^2+1^2+2^2}=n.$$
By (\ref{4.1}), we immediately see that $x+y+z+2w=s$. \qed

\begin{lemma}\label{Lem4.2} Let $x,y,z\in\Z$ with $2x^2+y^2+z^2\eq0\pmod 7$. Then $y^2\eq x^2\pmod 7$ or $z^2\eq x^2\pmod 7$.
\end{lemma}
\Proof. If $7\mid x$, then $y^2\eq-z^2\pmod 7$ and hence $y\eq z\eq0\eq x\pmod 7$ since $-1$ is a quadratic nonresidue modulo $7$.
In the case $7\nmid x$, there are $c,d\in\{0,1,2,4\}$ for which $y^2\eq cx^2\pmod 7$, $z^2\eq dx^2\pmod 7$ and $c+d\eq 5\pmod 7$,
thus $\{c,d\}=\{1,4\}$ and hence $x^2$ is congruent to $y^2$ or $z^2$ modulo $7$.
This ends the proof. \qed

\medskip
\noindent{\it Proof of Theorem} \ref{Th1.7}. (i) If $n\not\in E(1,1,1)$, then  for some $x,y,z,w\in\Z$ with $w=0$
we have $n=x^2+y^2+z^2+w^2$ and $w(x+y+z+2w)=0$.

Now assume that $n\in E(1,1,1)$. By (\ref{2.1}), $n=4^k(8l+7)$ for some $k,l\in\N$. As
$7n=4^k(8(7l+6)+1)\not\in E(1,1,1)$, we can write $7n=s^2+t^2+u^2+v^2$ with $s=0$ and $t,u,v\in\Z$.
As $t^2+u^2+v^2=7n-s^2\eq0\pmod 7$, we have $t^2\eq(2v)^2=4v^2\pmod 7$ and $u^2\eq(3v)^2\eq2v^2\pmod 7$, or $u^2\eq(2v)^2=4v^2\pmod 7$ and $t^2\eq(3v)^2\eq2v^2\pmod 7$.
Without any loss of generality, we simply assume that $t\eq2v\pmod7$ and $u\eq3v\pmod 7$. Define $x,y,z,w$ by (\ref{4.1}).
It is easy to see that $x,y,z,w\in\Z$.
In view of Lemma \ref{Lem4.1}, $x^2+y^2+z^2+w^2=n$ and $x+y+z+2w=s=0.$
Thus $w(x+y+z+2w)=0$ as desired.

 (ii) We can easily verify the desired result for all $n=0,\ldots,\lambda^214^{2m}-1$ via a computer.

Now let $n\gs \lambda^214^{2m}$ and assume the desired result for smaller values of $n$. If $4^m\mid n$, then by the induction hypothesis we have
$$\f n{4^m}=x^2+y^2+z^2+w^2\ \t{with}\ x,y,z,w\in\Z\ \t{and}\ x+y+z+2w\in\{\lambda r^m:\ r\in\N\},$$
and hence $n=(2^mx)^2+(2^my)^2+(2^mz)^2+(2^mw)^2$
with $$ 2^mx+2^my+2^mz+2\times 2^mw\in \{\lambda(2r)^m:\ r\in\N\}.$$

Now suppose that $4^m\nmid n$. By Lemma \ref{Lem3.5} we can write $7n$ as $(\lambda r^m)^2+t^2+u^2+v^2$ with $r\in\{0,7,14\}$ and $t,u,v\in\Z$.
Note that $s=\lambda r^m\eq0\pmod 7$. As in the proof of Theorem \ref{Th1.7}(i),
without any loss of generality we may assume that $t\eq2v\pmod7$ and $u\eq3v\pmod 7$. Define $x,y,z,w$ by (\ref{4.1}). Then $x,y,z,w\in\Z$.
By Lemma \ref{4.1}, $x^2+y^2+z^2+w^2=n$ and $x+y+z+2w=s=\lambda r^m.$

(iii) If $n\not\in E(1,1,1)$, then  for some $x,y,z,w\in\Z$ with $w=0$
we have $n=x^2+y^2+z^2+w^2$ and $w(x+2y+3z)=0$.

Now assume that $n\in E(1,1,1)$. By (\ref{2.1}), $n=4^k(8l+7)$ for some $k,l\in\N$.
As $7n=4^k(8(7l+6)+1)\not\in E(1,1,2)$ by (\ref{3.6}), we can write $7n=s^2+t^2+u^2+v^2$ with $s,t,u,v\in\Z$ with $s=t$.
As $u^2+v^2=7n-2s^2\eq5s^2\pmod 7$, without loss of generality we may assume that $u\eq2s\pmod 7$ and $v\eq-s\pmod 7$.
Clearly, those $x,y,z,w$ defined by (\ref{4.1}) are integers, and $x^2+y^2+z^2+w^2=n$ by Lemma \ref{Lem4.1}.
Note that $w+2x+3z=s-t=0$. Hence $n=w^2+x^2+z^2+y^2$ with $y(w+2x+3z)=0$.

(iv) For $n=0,1,\ldots,\lambda^2 14^{2m-1}-1$ we can verify the desired result via a computer.

Now fix an integer $n\gs \lambda^214^{2m-1}$ and assume that the desired result holds for smaller values of $n$.

If $4^m\mid n$, then by the induction hypothesis we can write $n/4^m=x^2+y^2+z^2+w^2$ $(x,y,z,w\in\Z)$
with $x+2y+3z=\lambda r^m$ for some $r\in\N$, hence
$$n=(2^mx)^2+(2^my)^2+(2^mz)^2+(2^mw)^2\ \ \t{with}\ 2^mx+2(2^my)+3(2^mz)=\lambda(2r)^m.$$

Below we suppose that $4^m\nmid n$. Applying Lemma \ref{Lem3.5} with $q=7$, we can write $14n=(\lambda r^m)^2+a^2+b^2+c^2$ with $r,a,b,c\in\N$.
Clearly, two of $a,b,c$ have the same parity. Without any loss of generality, we may assume that $b\eq c\pmod 2$ and $a\eq \lambda r^m\pmod 2$.
Set
$$s=\f{a+\lambda r^m}2,\ t=\f{a-\lambda r^m}2,\ u =\f{b+c}2,\ v=\f{b-c}2.$$
Then $s-t=\lambda r^m\eq0\pmod 7$. Note that
$$7n=\f{(\lambda r^m)^2+a^2}2+\f{b^2+c^2}2=s^2+t^2+u^2+v^2.$$
As $2t^2+u^2+v^2\eq0\pmod 7$, by Lemma \ref{Lem4.2}, one of $u^2$ and $v^2$ is congruent to $t^2$ modulo $7$, and another is congruent to $(2t)^2$ modulo 7.
Without any loss of generality, we may assume that $v\eq-t\pmod 7$ and $u\eq2t\pmod 7$. Now we define $x,y,z,w$ by (\ref{4.1}). It is easy to see that
$x,y,z,w\in\Z$. By Lemma \ref{Lem4.1}, $x^2+y^2+z^2+w^2=n.$
Note that $w+2x+3z=s-t=\lambda r^m$.

In view of the above, we have completed the proof of Theorem \ref{Th1.7}. \qed

\medskip
\noindent{\it Proof of Corollary} \ref{Cor1.2}. If $n\not\in E(1,1,1)$, then there are $w,x,y,z\in\Z$ with $w\not=0$ and $z=0$ such that $w^2+x^2+y^2+z^2=n$,
hence $x^2+24yz=x^2$ and $4x^2y^2+9y^2z^2+36z^2x^2=(2xy)^2$ are both squares.

Now assume that $n\in E(1,1,1)$. By Theorem \ref{Th1.7}(iii), $n$ can be written as $x^2+y^2+z^2+w^2\ (x,y,z,w\in\Z)$ with $w(x+2y+3z)=0$.
Since $n\in E(1,1,1)$, we have $w\not=0$ and $x+2y+3z=0$. Thus $x^2+y^2+(-z)^2+w^2=n$ with
$$x^2+24y(-z)=(2y+3z)^2-24yz=(2y-3z)^2.$$
Also, $n=y^2+x^2+z^2+w^2$ with
\begin{align*} 4y^2x^2+9x^2z^2+36z^2y^2
=&4x^2y^2+(x^2+4y^2)(3z)^2
\\=&(2xy)^2+(x^2+4y^2)(x^2+4xy+4y^2)
\\=&(x^2+2xy+4y^2)^2.
\end{align*}
The proof of Corollary \ref{Cor1.2} is now complete. \qed

\medskip
\noindent{\it Proof of Theorem} \ref{Th1.8}. (i) By Theorem \ref{Th1.7}(ii), there are $s,t,u,v\in\Z$ such that $5n=s^2+t^2+u^2+v^2$ and $s+t+u+2v=\lambda r^m$ for some $r\in\N$.
Note that a square is congruent to
one of $0,1,-1$ modulo $5$. As $s^2+t^2+u^2+v^2\eq0\pmod 5$, one of $s^2,t^2,u^2$ is congruent to $-v^2$ modulo $5$.
Without loss of generality, we assume that $u^2\eq-v^2\eq(2v)^2\pmod 5$ and hence $u+2v$ or $u-2v$ is divisible by $5$.
Since $(2s)^2\eq-s^2\eq t^2\pmod 5$, if $2s\not\eq t\pmod5$ then $t\eq -2s\pmod 5$ and hence $2t\eq-4s\eq s\pmod 5$.
Without any loss of generality, we simply assume that $2s\eq t\pmod 5$ and hence $2t\eq-s\pmod 5$. Define
$$x=\f{2s-t}5,\ y=\f{s+2t}5,\ z=\f{u+2v}5,\ w=\f{2u-v}5.$$
Then $x,y,5z,5w$ are all integers and
\begin{align*} x^2+y^2+z^2+w^2=&\f{(2s-t)^2+(s+2t)^2}{25}+\f{(u+2v)^2+(2u-v)^2}{25}
\\=&\f{s^2+t^2+u^2+v^2}5=n.\end{align*}
Note that
$$x+3y+5z=\f{(2s-t)+3(s+2t)}5+u+2v=s+t+u+2v=\lambda r^m.$$

(ii) By Lemma \ref{Lem3.3}, we can write $7n$ as $s^2+t^2+u^2+v^2$ $(s,t,u,v\in\Z)$
with $2s+t=\lambda r^m$ for some $r\in\Z$. Define $x,y,z,w$ as in (4.1). Then $7x,7y,7z,7w\in\Z$ and
$x^2+y^2+z^2+w^2=n$ by Lemma \ref{Lem4.1}.
Note that
$$x+3y+5w=2s+t=\lambda r^m.$$

The proof of Theorem \ref{Th1.8} is now complete. \qed
\medskip

\subsection*{Acknowledgements}
We would like to thank the referee for helpful comments.
This research was supported by the National Natural Science Foundation of China (grant no. 11571162)
and the NSFC-RFBR Cooperation and Exchange Program (grant 11811530072).

\end{document}